\documentclass{amsart}

\usepackage{latexsym}
\usepackage{amssymb}

\newtheorem{theorem}{Theorem}[section]
\newtheorem{lemma}[theorem]{Lemma}
\newtheorem{proposition}[theorem]{Proposition}

\theoremstyle{definition}
\newtheorem{definition}[theorem]{Definition}
\newtheorem{example}[theorem]{Example}

\theoremstyle{remark}

\numberwithin{equation}{section}

\begin{document}

\title[Proper cocycles and weak amenability]{Proper cocycles and weak forms of amenability}

\author{Paul Jolissaint}
\curraddr{\textsc{Universit\'e de Neuch\^atel,
       Institut de Math\'ematiques,      
       E.-Argand 11,
       CH-2000 Neuch\^atel, Switzerland}}
\email{paul.jolissaint@unine.ch}


\date{\today}



\begin{abstract}
Let $G$ and $H$ be locally compact, second countable groups. Assume that $G$ acts in a measure class preserving way on a standard space $(X,\mu)$ such that $L^\infty(X,\mu)$ has an invariant mean and that there is a Borel cocycle $\alpha:G\times X\rightarrow H$ which is proper in the sense of \cite{Jolcoc} and \cite{Knudby}.
We show that if $H$ has one of the three properties: Haagerup property (a-T-menability), weak amenability or weak Haagerup property, then so does $G$. In particular, we show that if $\Gamma$ and $\Delta$ are measure equivalent discrete groups in the sense of Gromov, then such cocycles exist and 
$\Gamma$ and $\Delta$ share the same weak amenability properties above.
\end{abstract}

\maketitle

\section{Introduction}

Let $G$ and $H$ be locally compact, second countable groups. If $H$ is a closed subgroup of $G$, it inherits weak amenability properties of $G$. Conversely, if in addition the homogeneous space $G/H$ is amenable in the sense that $L^\infty(G/H)$ has a $G$-invariant state, then $G$ in turn inherits weak amenability properties of $H$: see for instance \cite{AD}, \cite{ccjjv}, \cite{Jolcoc} and the recent article \cite{Knudby}. 

One way of proving such results is to use a natural cocycle associated to some regular Borel cross-section $\gamma:G/H\rightarrow G$ of the canonical projection $p:G\rightarrow G/H$. More precisely, we define $\alpha:G\times G/H\rightarrow H$ by
\[
\alpha(g,x)=\gamma(gx)^{-1}g\gamma(x)\quad ((g,x)\in G\times G/H).
\]
Since $g\gamma(x)H=\gamma(gx)H$ for all $g\in G$ and $x\in G/H$, we see that $\alpha(g,x)\in H$ as well. 
\par\vspace{3mm}
More generally, 
assume that $G$ acts in a measure class preserving way on some standard probability space $(X,\mu)$. Then a \textbf{Borel cocycle} is a Borel map
$\alpha:G\times X\rightarrow H$ such that for all $g,h\in G$, one has
\[
\alpha(gh,x)=\alpha(g,hx)\alpha(h,x)
\]
for $\mu$-almost every $x\in X$. Notice that $H$ is not necessarily a closed subgroup of $G$.

\par\vspace{3mm}
The aim of the present note is to prove that weak amenability properties of $H$ are inherited by $G$ when $\alpha$ is a cocycle that satisfies a properness condition in the sense of Definition \ref{1.3} below, and when $L^\infty(X)$ has a $G$-invariant state. In order to state precisely our main result, we need to recall some definitions from \cite{ccjjv}, \cite{CH}, \cite{Knudby}, \cite{Jolcoc} and \cite{Jolinv}. 

\par\vspace{3mm}
We assume throughout the article that our groups are locally compact and second countable. The definitions of the algebras $A(G)$ and $B_2(G)$ are reminded in the next section.

\begin{definition}
Let $G$ be a locally compact, second countable group.
\begin{enumerate}
\item [(1)] (\cite{ccjjv}) We say that $G$ has the \textbf{Haagerup property} if there exists a sequence $(u_n)_{n\geq 1}$ of normalized, positive definite functions on $G$, such that $u_n\in C_0(G)$ for every $n$, and that $u_n\to 1$ uniformly on compact subsets of $G$. ($C_0(G)$ denotes the vector space of all continuous functions on $G$ that tend to $0$ at infinity.)
\item [(2)] (\cite{CH}) The \textbf{Cowling-Haagerup constant} $\Lambda_{\mathrm{WA}}(G)$ is the infimum of all numbers $C\geq 1$ for which there exists a sequence $(u_n)_{n\geq 1}$ in the Fourier algebra $A(G)$ satisfying: $\Vert u_n\Vert_{B_2}\leq C$ for every $n$ and $u_n\to 1$ uniformly on compact sets, where $B_2(G)$ denotes the Herz-Schur multiplier algebra of $G$ (see \cite{Knudby}). Moreover, $G$ is said to be \textbf{weakly amenable} if $\Lambda_{\mathrm{WA}}(G)<\infty$.
\item [(3)] (\cite{Knudby}) The \textbf{weak Haagerup constant} $\Lambda_{\mathrm{WH}}(G)$ is the infimum of all numbers $C\geq 1$ for which there exists a sequence $(u_n)_{n\geq 1}\subset B_2(G)\cap C_0(G)$ such that $\Vert u_n\Vert_{B_2}\leq C$ for every $n$ and $u_n\to 1$ uniformly on compact sets. Moreover, $G$ is said to have the \textbf{weak Haagerup property} if $\Lambda_{\mathrm{WH}}(G)<\infty$.
\end{enumerate}
\end{definition}

As $A(G)\subset B_2(G)\cap C_0(G)$, one always has $\Lambda_{\mathrm{WH}}(G)\leq \Lambda_{\mathrm{WA}}(G)$, and a weakly amenable group has the weak Haagerup property. Similarly, as normalized, positive definite functions are Herz-Schur multipliers of norm one, if $G$ has the Haagerup property then it has the weak Haagerup property and $\Lambda_{\mathrm{WH}}(G)=1$. See \cite{Knudby} for a discussion of these properties.

The following definition generalizes the notion of \textbf{co-F\o lner groups} as in \cite{ccjjv}.

\begin{definition}\label{ampair}
(\cite{Gre}, \cite{Jolinv}, \cite{Zimmer}) 
Let $G$ be a locally compact, second countable group that acts in a measure class preserving way on the standard space $(X,\mu)$. Then we say that $(G,X)$ is an \textbf{amenable pair} if $L^\infty(X,\mu)$ has a $G$-invariant state.
\end{definition}

Finally, let $\alpha:G\times X\rightarrow H$ be a Borel cocycle. The following is essentially taken from \cite{Jolcoc}; see also \cite{Knudby}. We need first to fix some notation: let $A\subset X$ be a Borel set and let $L$ be a compact subset of $H$; we denote by $K(A,L)$ the set of all elements $g\in G$ for which $\mu(X_A(g))>0$, where 
\[
X_A(g)=\{x\in A\cap g^{-1}A: \alpha(g,x)\in L\}.
\]

As mentioned above, it is motivated by the case where $H$ is a closed subgroup of $G$ that we present in detail now. Choose a regular Borel cross-section $\gamma:G/H\rightarrow G$ for the canonical projection $p:G\rightarrow G/H$, \textit{i.e.} $\gamma$ is a Borel map such that $p(\gamma(x))=x$ for every $x\in G/H$ and such that, for every compact set $K\subset G$, the set $\gamma(G/H)\cap p^{-1}(p(K))$ is precompact in $G$. Moreover, for every compact set $C\subset G/H$, there is a compact set $C_1\subset G$ such that $C=p(C_1)$. Then it is straightforward to see that $\gamma(C)$ is precompact for every compact set $C\subset G/H$. 

Recall that the associated cocycle $\alpha$ is defined by
\[
\alpha(g,x)=\gamma(gx)^{-1}g\gamma(x)\quad (g\in G, x\in G/H).
\]
Then one has:
\begin{enumerate}
\item [(i)] If $K\subset G$ and $A\subset G/H$ are compact sets, then $\alpha(K\times A)$ is precompact since it is contained in  $\gamma(KA)^{-1}K\gamma(A)$.
\item [(ii)] If $A\subset G/H$ and $L\subset H$ are compact, it is easy to see that, if $g\in G$ is such that $X_A(g)\not=\emptyset$, then $g\in \gamma(A)L\gamma(A)^{-1}$. In particular $K(A,L)$ is precompact.
\end{enumerate}
As $G$ and $G/H$ are in particular $\sigma$-compact, we see that the cocycle $\alpha$ satisfies the next definition.

\begin{definition}\label{1.3}
The cocycle $\alpha:G\times X\rightarrow H$ is \textbf{proper} if it satisfies the following two conditions:
\begin{enumerate}
\item [(i)] for every Borel set $A\subset X$, for every $\varepsilon>0$ and for every compact set $K\subset G$, there exists a Borel set $A_{\varepsilon,K}\subset A$ such that $\mu(A\smallsetminus A_{\varepsilon,K})<\varepsilon$ and 
$\alpha(K\times A_{\varepsilon,K})$ is precompact.
\item [(ii)] For every compact set $L\subset H$ and every $\varepsilon>0$, there exists a Borel set $A\subset X$ such that $\mu(X\smallsetminus A)\leq\varepsilon$ and $K(A,L)$ is precompact.
\end{enumerate}
\end{definition}

We observe first that properness of $\alpha$ is independent of the chosen probability measure.

\begin{lemma}
Let $\alpha:G\times X\rightarrow H$ be a proper cocycle with respect to the probability measure $\mu$ on $X$, and let $\nu$ be an equivalent probability measure on $X$. Then $\alpha$ is proper with respect to $\nu$.
\end{lemma}
\textit{Proof.} This follows immediately from Theorem 6.11 of \cite{Rudin}: as $\nu$ is equivalent to $\mu$, it follows in particular that for every $\varepsilon>0$, there exists $\theta>0$ such that, if $B\subset X$ is Borel and if $\mu(B)\leq\theta$, then $\nu(B)\leq \varepsilon$.
\hfill $\square$

\par\vspace{3mm}
Before giving examples of such cocycles, let us state our first main result:

\begin{theorem}\label{thmproper}
Let $G$ and $H$ be locally compact, second countable groups, let $G$ act on some probability space $(X,\mu)$ so that $(G,X)$ is an amenable pair, and let $\alpha:G\times X\rightarrow H$ be a proper cocycle. 
\begin{enumerate}
\item [(a)] If $H$ has the Haagerup property, then so does $G$.
\item [(b)] If $H$ is weakly amenable group then so is $G$, and 
\[
\Lambda_{\mathrm{WA}}(G)\leq \Lambda_{\mathrm{WA}}(H).
\]
\item [(c)] If $H$ has the weak Haagerup property, then so does $G$, and
\[
\Lambda_{\mathrm{WH}}(G)\leq \Lambda_{\mathrm{WH}}(H).
\]
\end{enumerate}
\end{theorem}

As we will see, the proofs of the three statements rely on the same techniques, and they will be given in the next section. 
\par\vspace{3mm}
As promised, here are examples of proper cocycles.

\begin{example}
Every cocycle $\alpha$ described below is proper.
\begin{enumerate}
\item [(1)] Let $G$ and $H$ be locally compact, second countable groups and assume that $\sigma: G\rightarrow H$ is a continuous homomorphism with compact kernel, and let $(X,\mu)$ be an arbitrary standard probability $G$-space. Define $\alpha:G\times X\rightarrow H$ by
\[
\alpha(g,x)=\sigma(g)\quad ((g,x)\in G\times X).
\]
\item [(2)] More generally, let $G$ and $H$ be as in Example (1), let $G_0$ be a closed subgroup of $G$ and assume that $\sigma:G_0\rightarrow H$ is a continuous homomorphism with compact kernel. Choose a regular Borel cross-section $\gamma:G/G_0\rightarrow G$ for the canonical projection and a quasi-invariant probability measure $\nu$ on $G/G_0$. If $(Y,\mu)$ is an arbitrary standard probability $G$-space, equip $(G/G_0\times Y,\nu\times \mu)$ with the product action and define $\alpha:G\times (G/G_0\times Y)\rightarrow H$ by
\[
\alpha(g,(x,y))=\sigma(\gamma(gx)^{-1}g\gamma(x)).
\]
This case generalizes the situation where $H=G_0$ is a closed subgroup of $G$.
\item [(3)] Let $\pi:P\rightarrow B$ be a (metrizable, locally compact) topological principal fiber bundle with countable structure group $H$, and assume that a locally compact, second countable group $G$ acts continuously on $P$ so that the $G$-action commutes with the $H$-action. Using a bounded measurable cross-section for $\pi$, we get a trivialization of $P\cong H\times B$ which preserves precompact subsets. Then the corresponding action of $G$ on $H\times B$ is given by
\[
g\cdot(h,x)=(\alpha(g,x)h,gx)
\]
where $\alpha$ is a Borel cocycle. If the action of $G$ on $P$ is proper, then $\alpha$ is a proper cocycle.
\end{enumerate}
\end{example}

As will be proved in the last section, an important fourth family of pairs of groups that give rise to proper cocycles is the family of pairs of countable, discrete groups $\Gamma$ and $\Delta$ that satisfy Gromov's notion of measure equivalence.
We recall the latter from \cite{Furman1}:

\begin{definition}
We say that $\Gamma$ and $\Delta$ are \textbf{measure equivalent} if there exist commuting, measure-preserving, free actions of $\Gamma$ and $\Delta$ on some infinite-measure standard space $(\Sigma,\sigma)$, such that $\Gamma$ and $\Delta$ both admit fundamental domains with finite measure. (For convenience, we denote the action of $\Gamma$ on the left and the action of $\Delta$ on the right.)
\end{definition}

Measure equivalence is a weak form of orbit equivalence: see \cite{Furman2}, Section 3. Let us generalize it slightly as follows.

\begin{definition}
We say that $\Gamma$ and $\Delta$ are \textbf{amenably measure equivalent} if there exist commuting, measure-preserving, free actions of $\Gamma$ and $\Delta$ on some infinite-measure standard space $(\Sigma,\sigma)$, such that the pairs $(\Gamma,\Sigma/\Delta)$ and $(\Delta,\Gamma\backslash\Sigma)$ are both amenable pairs.
\end{definition}

\begin{example}
Let $\Gamma$ and $\Delta$ be discrete subgroups of the same locally compact, second countable, unimodular group $G$ such that the homogeneous spaces $G/\Gamma$ and $G/\Delta$ both have $G$-invariant means, \textit{i.e.} the pairs $(G,G/\Gamma)$ and $(G,G/\Delta)$ are amenable.
Then $\Gamma$ and $\Delta$ are amenably measure equivalent groups.
\end{example}

As claimed above, amenably measure equivalent groups give rise to proper cocycles, so that Theorem \ref{thmproper} will be used to prove our last result.

\begin{theorem}\label{thmAME}
Let $\Gamma$ and $\Delta$ be amenably measure equivalent groups. 
\begin{enumerate}
\item [(1)] If one of them has the Haagerup property, then the other one has the same property.
\item [(2)] The following equalities hold:
\[
\Lambda_{\mathrm{WA}}(\Gamma)=\Lambda_{\mathrm{WA}}(\Delta)\quad\textrm{and}\quad
\Lambda_{\mathrm{WH}}(\Gamma)=\Lambda_{\mathrm{WH}}(\Delta).
\]
\end{enumerate}
\end{theorem}

Let us discuss instances where Theorem \ref{thmAME} applies in the context of semidirect products. Let $A$ be an amenable, countable group. If it acts on some group $\Gamma$, then it is easy to see that $\Gamma$ and $\Gamma\rtimes A$ are amenably measure equivalent. 

The fact that $\Gamma$ and $\Gamma\rtimes A$ have simultaneously the Haagerup property and that their constants $\Lambda_{\mathrm{WA}}$ and $\Lambda_{\mathrm{WH}}$ coincide is already known: it follows respectively from \cite{ccjjv}, \cite{Jolcoc} and \cite{Knudby}.
\par\vspace{3mm}
Next, let us assume that the amenable group $A$ is on the other side:
let
$\theta:\Delta\rightarrow \mathrm{Aut}(A)$ be an action of the countable
group $\Delta$ on $A$, and consider the following action of the semidirect
product group $\Gamma:=A\rtimes_{\theta}\Delta$ on $A$: 
\[
(a,\delta)\cdot b=a\theta_{\delta}(b)
\] 
for all $a,b\in A$ and all
$\delta\in\Delta$. Then we have:

\begin{proposition}\label{semidirect}
Let $A$, $\Delta$ and $\theta$ be as above and assume that there
exists a sequence $(F_k)_{k\geq 1}$ of finite subsets of $A$ such
that, for all $(a,\delta)\in \Gamma$ $$
\frac{|a\theta_{\delta}(F_k)\bigtriangleup F_k|}{|F_k|}\to 0$$ as
$k\to\infty$. In other words, the pair $(\Gamma,A)$ is amenable in
the sense of Definition \ref{ampair}. Then $\Gamma$ and $\Delta$ are
amenably measure equivalent.
\end{proposition}
\textit{Proof}. Take $\Sigma=A\times\Delta$ with counting measure.
The action (on the left) of $\Gamma=A\rtimes\Delta$ is defined by
multiplication in the semidirect product, namely 
\[
(a,\delta)\cdot(b,\delta')=(a\theta_{\delta}(b),\delta\delta')
\]
and the action of $\Delta$ (on the right) is given by
$(a,\delta)\cdot\delta'=(a,\delta\delta')$. Then
$\Gamma\backslash\Sigma$ is the one-point space and
$\Sigma/\Delta$ is isomorphic to $A$ as a $\Gamma$-space. The
amenability of the pair $(\Gamma,A)$ implies that $\Gamma$ and $\Delta$ are
amenably measure equivalent. \hfill $\square$

\par\vspace{3mm}
Let us illustrate Proposition \ref{semidirect} by an example.

\begin{example}
Let $1\rightarrow Z\rightarrow A\rightarrow Q\rightarrow 1$
be a central extension of some amenable group $Q$. Denote by
$s:Q\rightarrow A$ a section of the canonical projection
with $s(1)=1$. 

Let also $h:\Delta\rightarrow Q$ be a
homomorphism of some group $\Delta$. With these data, define an
action $\theta$ of $\Delta$ on $A$ by
\[
\theta_{\delta}=\mathrm{Ad}(s\circ h(\delta))\quad\forall\delta\in\Delta.
\]
As each automorphism $\theta_{\delta}$ is inner,
$(A,\Delta,\theta)$ fulfills the conditions of Proposition \ref{semidirect}:
indeed, since $A$ is amenable, the direct product group $A\times A$ is as well, and the action $(a,b)\cdot x=axb^{-1}$
of $A\times A$ on $A$ admits a sequence $(F_k)_{k\geq 1}$ of finite subsets
of $A$ such that 
\[
\frac{|(a,b)\cdot F_k\bigtriangleup
F_k|}{|F_k|}\to 0
\] 
as $k\to\infty$ for all $a,b\in A$.

Remark also that $s\circ h$ is not a homomorphism in general,
thus $\Gamma:=A\rtimes_{\theta}\Delta$ is not isomorphic to the direct
product $A\times \Delta$.
\hfill $\square$
\end{example}

\par\vspace{2mm}
\textit{Acknowledgements.} I am very grateful to Tadeusz Januszkiewicz and the referee for their careful reading of the manuscript and their valuable comments.


\section{Proof of Theorem \ref{thmproper}}

Let $G$ and $H$ be locally compact groups; we assume that they are second countable even if definitions below make sense for arbitrary locally compact groups. 

The \textbf{Fourier-Stieltjes algebra} of $G$ is the set of all coefficient functions of unitary representations of $G$, thus, for every $u\in B(G)$ there exists a unitary representation $(\pi,\mathcal H)$ of $G$ and two vectors $\xi,\eta\in\mathcal H$ such that $u(g)=\langle\pi(g)\xi|\eta\rangle$ for every $g\in G$. It is a Banach algebra with respect to the norm
\[
\Vert u\Vert_B=\inf\Vert\xi\Vert\Vert\eta\Vert
\]
where the infimum is taken over all representations of $u$ as above.

The \textbf{Fourier algebra} of $G$ is the set of all coefficient functions associated to the left regular representation $\lambda$ of $G$ (which acts on $L^2(G)$).
It is the norm closure of the algebra of compactly supported continuous functions $C_c(G)\cap B(G)$ in the  algebra $B(G)$.

A \textbf{Herz-Schur multiplier} of $G$ is a continous function $u:G\rightarrow\mathbb C$ for which there exists a separable Hilbert space $\mathcal H$ and two bounded, continuous functions $\xi,\eta:G\rightarrow\mathcal H$ such that
\[
u(h^{-1}g)=\langle \xi(g)|\eta(h)\rangle\quad (g,h\in G).
\]
It turns out that the set $B_2(G)$ of all Herz-Schur multipliers on $G$ is a Banach algebra with respect to the pointwise product and to the norm
\[
\Vert u\Vert_{B_2}=\inf\Vert \xi\Vert_\infty\Vert\eta\Vert_\infty
\]
where the infimum is taken over all representations of $u$ as above.

\par\vspace{3mm}

Assume from now on that $G$ acts in a measure class preserving way on a standard probability space $(X,\mu)$ and that $\alpha:G\times X\rightarrow H$ is a (not necessarily proper) Borel cocycle.

We denote by $(g,x)\mapsto \chi(g,x)$ the Radon-Nikodym derivative related to the action of $G$ on $X$ and characterized by
\[
\int\limits_X f(gx)\chi(g,x)d\mu(x)=\int\limits_X f(x)d\mu(x)\quad\forall f\in L^1(X,\mu).
\]
It satisfies the cocycle relation:
\[
\chi(gh,x)=\chi(g,hx)\chi(h,x)
\]
for all $g,h\in G$ and $\mu$-a.e. $x\in X$. Taking $f=1_X$, we have: $\int\limits_X\chi(g,x)d\mu(x)=1$ for every $g\in G$.
\par\vspace{3mm}
For future use, let us observe that for every $g\in G$ and every Borel set $B\subset X$, one has, by Cauchy-Schwarz Inequality:
\[
\int\limits_B\sqrt{\chi(g,x)}d\mu(x)\leq \mu(B)^{1/2}\left(\int\limits_X\chi(g,x)d\mu(x)\right)^{1/2}=\mu(B)^{1/2}
\]
and
\[
\int\limits_{g^{-1}B}\sqrt{\chi(g,x)}d\mu(x)\leq \mu(X)^{1/2}\left(\int\limits_X 1_B(gx)\chi(g,x)d\mu(x)\right)^{1/2}=\mu(B)^{1/2}.
\]

\par\vspace{3mm}

The proof of Theorem \ref{thmproper} relies on two auxiliary results.

\begin{lemma}\label{2.1}
Let $u\in L^\infty(H)$. Define $\hat u:G\rightarrow \mathbb C$ by
\[
\hat u(g)=\int\limits_X u(\alpha(g,x))\sqrt{\chi(g,x)}d\mu(x)\quad\forall g\in G.
\]
\begin{enumerate}
\item [(a)] If $u\in B_2(H)$ is a Herz-Schur multiplier on $H$, then $\hat u\in B_2(G)$ and $\Vert\hat u\Vert_{B_2}\leq \Vert u\Vert_{B_2}$. If furthermore $u$ is positive definite, so is $\hat u$.
\item [(b)] If $\alpha$ is proper and if $u\in C_0(H)$ then $\hat u\in C_0(G)$. In particular, if $u\in B_2(H)\cap C_0(H)$, then $\hat u\in B_2(G)\cap C_0(G)$.
\item [(c)] If $\alpha$ is proper and if $u\in A(H)$, then $\hat u\in A(G)$.
\end{enumerate}
\end{lemma}
\textit{Proof.} (a) There exist a separable Hilbert space $\mathcal H$ and bounded, continuous functions $\xi,\eta:H\rightarrow\mathcal H$ such that
\begin{enumerate}
\item [(1)] $u(t^{-1}s)=\langle\xi(s)|\eta(t)\rangle$ for all $s,t\in H$;
\item [(2)] $\Vert u\Vert_{B_2}\leq \Vert \xi\Vert_\infty \Vert \eta\Vert_\infty$.
\end{enumerate}
Define $\hat\xi,\hat\eta:G\rightarrow L^2(X,\mu,\mathcal H)$ by
\[
\hat\xi(g)(x)=\xi(\alpha(g^{-1},x)^{-1})\sqrt{\chi(g^{-1},x)}
\]
and
\[
\hat\eta(g)(x)=\eta(\alpha(g^{-1},x)^{-1})\sqrt{\chi(g^{-1},x)}.
\]
One has for every $g\in G$:
\begin{eqnarray*}
\Vert\hat\xi(g)\Vert^2 &=& \int\limits_X\Vert\xi(\alpha(g^{-1},x))\Vert^2\chi(g^{-1},x)d\mu(x)\\
&\leq &
\Vert\xi\Vert_\infty^2\int\limits_X\chi(g^{-1},x)d\mu(x)=\Vert\xi\Vert_\infty^2.
\end{eqnarray*}
Similarly, $\Vert\hat\eta\Vert_\infty\leq\Vert\eta\Vert_\infty$. We are going to prove that, for all $g,h\in G$, one has:
\[
\hat u(h^{-1}g)=\langle\hat\xi(g)|\hat\eta(h)\rangle.
\]
It turns out that $\hat u$ is a continuous function on $G$ by Appendix C of \cite{Knudby} (even though $\hat\xi$ and $\hat\eta$ are not necessarily continuous).

Observe that the cocycle relation $\alpha(gh,x)=\alpha(g,hx)\alpha(h,x)$ for all $g,h\in G$ and for $\mu$-a.e $x\in X$ implies that
\[
\alpha(h^{-1}g,g^{-1}x)\alpha(g^{-1},x)=\alpha(h^{-1},x)
\] and similarly
\[
\chi(h^{-1}g,g^{-1}x)\chi(g^{-1},x)=\chi(h^{-1},x)
\]
for all $g,h\in G$ and $\mu$-a.e. $x\in X$. Fix $g,h\in G$. One has:
\begin{eqnarray*}
\langle\hat\xi(g)|\hat\eta(h)\rangle &=&
\int\limits_X\langle \xi(\alpha(g^{-1},x)^{-1})|\eta(\alpha(h^{-1},x)^{-1})\rangle\sqrt{\chi(g^{-1},x)\chi(h^{-1},x)}d\mu(x)\\
&=&
\int\limits_X u(\alpha(h^{-1},x)\alpha(g^{-1},x)^{-1})\sqrt{\chi(g^{-1},x)\chi(h^{-1},x)}d\mu(x)\\
&=&
\int\limits_X u(\alpha(h^{-1}g,g^{-1}x))\sqrt{\chi(g^{-1},x)\chi(h^{-1},x)}d\mu(x)\\
&=&
\int\limits_X u(\alpha(h^{-1}g,g^{-1}x))\chi(g^{-1},x)\sqrt{\chi(h^{-1}g,g^{-1}x)}d\mu(x)\\
&=&
\int\limits_X u(\alpha(h^{-1}g,x))\sqrt{\chi(h^{-1}g,x)}d\mu(x)=\hat u(h^{-1}g).
\end{eqnarray*}
If furthermore $u$ is positive definite, then one can take $\eta=\xi$, and it is straightforward to see that $\hat u$ is positive definite on $G$ as well. (In fact, let $(\pi_u,\mathcal{H}_u,\xi_u)$ be the Gel'fand-Naimark-Segal triple associated to $u$. Then, as $u(h)=\langle \pi_u(h)\xi_u|\xi_u\rangle$ for every $h\in H$, we see that the function $h\mapsto \xi(h)=\pi_u(h)\xi_u$ works.)\\
(b) Let now $u\in C_0(H)$. We assume without loss of generality that $\Vert u\Vert_\infty\leq 1$. 

Fix $\varepsilon>0$. There exists a compact set $L\subset H$ such that
$|u(h)|\leq\frac{\varepsilon}{2}$ for all $h\notin L$.
Choose next a Borel set $A\subset X$ which satisfies: $\mu(X\smallsetminus A)\leq \frac{\varepsilon^2}{16}$ and for which $K=\overline{K(A,L)}$ is compact.

Fix $g\in G\smallsetminus K$. One has:
\begin{eqnarray*}
|\hat u(g)| &\leq&
\int\limits_{\{\alpha(g,x)\in L\}} |u(\alpha(g,x))|\sqrt{\chi(g,x)}d\mu(x)\\
& &
+\int\limits_{\{\alpha(g,x)\notin L\}}|u(\alpha(g,x))|\sqrt{\chi(g,x)}d\mu(x)\\
&\leq&
\int\limits_{X\smallsetminus(A\cap g^{-1}A)}|u(\alpha(g,x))|\sqrt{\chi(g,x)}d\mu(x)+\frac{\varepsilon}{2}\\
&\leq&
\int\limits_{X\smallsetminus A}\sqrt{\chi(g,x)}d\mu(x)+\int\limits_{g^{-1}(X\smallsetminus A)}\sqrt{\chi(g,x)}d\mu(x)+\frac{\varepsilon}{2}\\
&\leq&
2\mu(X\smallsetminus A)^{1/2}+\frac{\varepsilon}{2}\leq\varepsilon.
\end{eqnarray*}
(c) The assertion is similar to that contained in Proposition 2.8 of \cite{Jolcoc}, but we give a proof for the sake of completeness. We have to prove that $\hat u$ is in the norm closure of $C_c(G)\cap B(G)$ in the Fourier-Stieltjes algebra $B(G)$. Thus, let us assume that $u(h)=\langle \lambda(h)\xi|\eta\rangle$ with $\xi,\eta\in C_c(H)$. Then it is straighforward to check that 
\[
\hat u(g)=\langle\lambda_\alpha(g)1_X\otimes\xi|1_X\otimes\eta\rangle\quad (g\in G)
\]
where $\lambda_\alpha$ is the unitary representation of $G$ on $L^2(X,\mu,L^2(H))$ defined by
\[
(\lambda_\alpha(g)\zeta)(x)=\lambda(\alpha(g^{-1},x)^{-1})\zeta(g^{-1}x)\sqrt{\chi(g^{-1},x)}
\]
for $\zeta\in L^2(X,\mu,L^2(H))$. Hence $\hat u$ is a norm limit in $B(G)$ of functions of the form 
\[
\hat{u}_A(g)=\langle\lambda_\alpha(g)1_A\otimes\xi|1_A\otimes\eta\rangle
\]
with $A\subset X$ Borel. Denote by $L$ the support of $u$. For every Borel set $A\subset X$ such that $K(A,L)$ is precompact, we have
\[
\hat{u}_A(g)=\int\limits_{A\cap g^{-1}A}u(\alpha(g,x))\sqrt{\chi(g,x)}d\mu(x) \quad (g\in G).
\]
If $g\notin\overline{K(A,L)}$, the latter being compact in $G$, then $\mu(\{x\in A\cap g^{-1}A:\alpha(g,x)\in L\})=0$ by Definition \ref{1.3} and this implies that $\hat{u}_A(g)=0$.
\hfill $\square$

\par\vspace{3mm}
Assume now that $(G,X)$ is an amenable pair. Denote by $\beta$ the action
of $G$ on $L^1(X,\mu)$ given by
\[
\beta_g(f)(x)=f(g^{-1}x)\chi(g^{-1},x)
\]
Since the set of normal states is weak* dense in the set of all states of $L^\infty(X)$, the amenability of $(G,X)$ is equivalent to the existence of a sequence $(f_n)\subset L^1(X,\mu)$ such that
\begin{enumerate}
\item [(1)] $f_n\geq 0$ and $\Vert f_n\Vert_1=\int\limits_X f_n(x)d\mu(x)=1$ for every $n$;
\item [(2)] for every compact set $K\subset G$, 
\[
\sup_{g\in K}\Vert \beta_g(f_n)-f_n\Vert_1\to 0
\]
as $n\to\infty$.
\end{enumerate}

\begin{lemma}
If the pair $(G,X)$ is amenable, then there exists a sequence of probability measures $(\mu_n)$ on $X$ such that:
\begin{enumerate}
\item [(a)] $\mu_n$ is equivalent to $\mu$ for every $n$;
\item [(b)] for every compact set $K\subset G$,
\[
\sup_{g\in K}\int\limits_X |\sqrt{\chi_n(g,x)}-1|d\mu_n(x)\to 0
\]
as $n\to\infty$, where $\chi_n(g,\cdot)=dg_*^{-1}\mu_n/d\mu_n$ denotes the corresponding Radon-Nikodym derivative.
\end{enumerate}
\end{lemma}
\textit{Proof.} Let $(f_n)$ be a sequence in $L^1(X,\mu)$ as above. Adding the constant function $\frac{1}{n}1_X$ to $f_n$ and renormalizing if necessary, we assume that $(f_n)$ satisfies conditions (1) and (2) above, and that there exists a constant $c_n>0$ for every $n$ such that $f_n\geq c_n$ for every $n$. Define $\mu_n$ by
\[
\int\limits_Xf(x)d\mu_n(x)=\int\limits_X f(x)f_n(x)d\mu(x).
\]
Then $L^\infty(X,\mu_n)=L^\infty (X,\mu)$ isometrically, and we have for $f\in L^\infty(X,\mu)$, $g\in G$ and $n\geq 1$:
\begin{eqnarray*}
\left|\int\limits_X f(x)\lbrace \chi_n(g,x)-1\rbrace d\mu_n(x)\right| &=&
\left|\int\limits_X f(g^{-1}x)f_n(x)d\mu(x)-\int\limits_X f(x)f_n(x)d\mu(x)\right|\\
&=&
\left|\int\limits_X f(x)\lbrace f_n(gx)\chi(g,x)-f_n(x)\rbrace d\mu(x)\right|\\
&=&
\left|\int\limits_X f(x)\lbrace\beta_{g^{-1}}(f_n)(x)-f_n(x)\rbrace d\mu(x)\right|\\
&\leq &
\Vert f\Vert_\infty \Vert \beta_{g^{-1}}(f_n)-f_n\Vert_1.
\end{eqnarray*}
This implies that 
\[
\Vert\chi_n(g,\cdot)-1\Vert_{L^{1}(X,\mu_n)}=\int\limits_X|\chi_n(g,x)-1|d\mu_n(x)\leq \Vert\beta_{g^{-1}}(f_n)-f_n\Vert_{L^{1}(X,\mu)}
\]
for every $g\in G$ and every $n$. Let $K$ be a compact subset of $G$ and $\varepsilon>0$. There exists $N>0$ such that
\[
\sup_{g\in K}\Vert\beta_{g^{-1}}(f_n)-f_n\Vert_{L^{1}(X,\mu)}\leq\varepsilon
\]
for every $n\geq N$, so that we get for $g\in K$ and $n\geq N$:
\begin{eqnarray*}
\int\limits_X|\sqrt{\chi_n(g,x)}-1|d\mu_n(x) &\leq &
\int\limits_X |\sqrt{\chi_n(g,x)}-1||\sqrt{\chi_n(g,x)}+1|d\mu_n(x)\\
&=&
\int\limits_X|\chi_n(g,x)-1|d\mu_n(x)\leq\varepsilon.
\end{eqnarray*}
\hfill $\square$

\par\vspace{3mm}\noindent
\textit{Proof of Theorem \ref{thmproper}.} We are going to prove statement (c) of Theorem 1.4. The proofs of (a) and (b) are special cases that will be discussed briefly afterwards. 

Thus let us assume that $H$ has the weak Haagerup property. 

Fix $C>\Lambda_{\mathrm{WH}}(H)$, $K\subset G$ compact and $\varepsilon>0$. We are going to prove that there exists $\hat u\in B_2(G)\cap C_0(G)$ such that:
\begin{enumerate}
\item [(1)] $\Vert\hat u\Vert_{B_2}\leq C$;
\item [(2)] $\sup_{g\in K}|\hat u(g)-1|\leq\varepsilon$.
\end{enumerate}
Let $(\mu_n)$ be as in Lemma 2.2, and let $n$ be large enough in order that 
\[
\sup_{g\in K}\int\limits_X|\sqrt{\chi_n(g,x)}-1|d\mu_n(x)\leq \frac{\varepsilon}{4(C+1)}=:\varepsilon'.
\]
By condition (i) of Definition \ref{1.3} and Lemma 1.4, there exists a Borel set $A=A_{\varepsilon',K}$ such that $\alpha(K\times A)$ is precompact and 
\[
\mu_n(X\setminus A)\leq\frac{\varepsilon}{4(C+1)}.
\]
Let $L$ denote the closure of $\alpha(K\times A)$ in $H$ and choose $u\in B_2(H)\cap C_0(H)$ such that
\[
\Vert u\Vert_{B_2}\leq C\quad\textrm{and}\quad
\sup_{h\in L}|u(h)-1|\leq\frac{\varepsilon}{4}
\]
and put
\[
\hat u(g)=\int\limits_X u(\alpha(g,x))\sqrt{\chi_n(g,x)}d\mu_n(x)
\]
for $g\in G$. Lemma \ref{2.1} implies that $\hat u\in B_2(G)\cap C_0(G)$ and that $\Vert \hat u\Vert_{B_2}\leq C$.

We have for every $g\in K$:
\begin{eqnarray*}
|\hat u(g)-1| &\leq&
\int\limits_{A}|u(\alpha(g,x))\sqrt{\chi_n(g,x)}-1|d\mu_n(x)\\
& &
+\int\limits_{A^c}|u(\alpha(g,x))\sqrt{\chi_n(g,x)}-1|d\mu_n(x).
\end{eqnarray*}
Then
\begin{eqnarray*}
\int\limits_{A}|u(\alpha(g,x))\sqrt{\chi_n(g,x)}-1|d\mu_n(x)
& \leq &
\int\limits_{A}|u(\alpha(g,x))||\sqrt{\chi_n(g,x)}-1|d\mu_n(x)\\
&  &
+\int\limits_{A}|u(\alpha(g,x))-1|d\mu_n(x)\\
&\leq&
C\cdot \int\limits_X|\sqrt{\chi_n(g,x)}-1|d\mu_n(x)+\frac{\varepsilon}{4}\\
& \leq & 
\frac{\varepsilon}{2}.
\end{eqnarray*}
Next,
\[
\int\limits_{A^c}|u(\alpha(g,x))\sqrt{\chi_n(g,x)}-1|d\mu_n(x)
\]
\begin{eqnarray*}
& \leq &
\int\limits_{A^c} |u(\alpha(g,x))||\sqrt{\chi_n(g,x)}-1|d\mu_n(x)\\
& &
+\int\limits_{A^c} |u(\alpha(g,x))-1|d\mu_n(x)\\
&\leq&
C\cdot \int\limits_X |\sqrt{\chi_n(g,x)}-1|d\mu_n(x)\\
& &
+(C+1)\mu_n(X\smallsetminus X(K,L))\\
&\leq & \frac{\varepsilon}{2}.
\end{eqnarray*}
This ends the proof of statement (c).

If $H$ satisfies condition (a), given $K\subset G$ compact and $\varepsilon>0$,
the same construction as above from a positive definite, normalized function $u\in C_0(H)$ gives a positive definite function $\hat u\in C_0(G)$ that satisfies 
\[
\sup_{g\in K}|\hat u(g)-1|\leq\varepsilon.
\]
This proves that $G$ has the Haagerup property.

Finally, if $H$ satisfies condition (b), if $C>\Lambda_{\mathrm{WA}}(H)$, $K\subset G$ compact and $\varepsilon>0$ are given, choosing $u\in A(H)$ with $\Vert u\Vert_{B_2}\leq C$ as in the first part of the proof, we get $\hat u\in A(G)$ satisfying (1) and (2) above. This proves that $G$ is weakly amenable and that $\Lambda_{\mathrm{WA}}(G)\leq\Lambda_{\mathrm{WA}}(H)$.
\hfill $\square$

\section{Proof of theorem \ref{thmAME}}

Let $\Gamma$ and $\Delta$ be amenably measure equivalent groups as in Theorem \ref{thmAME} and let $(\Sigma,\sigma)$ be a standard infinite measure space on which $\Gamma$ and $\Delta$ act, the first one acting freely on the left and the second one freely on the right, both actions preserving the infinite measure $\sigma$, and such that the pairs $(\Gamma,\Sigma/\Delta)$ and $(\Delta,\Gamma\backslash\Sigma)$ are amenable.

We fix our notation and recall some needed facts from \cite{Furman1}: we choose Borel cross-sections from $\Gamma\backslash\Sigma$ and $\Sigma/\Delta$ to $\Sigma$, and we denote by $Y$ and $X$ their respective ranges so that
\[
\Sigma=\bigsqcup_{\delta\in\Delta}X\delta=
\bigsqcup_{\gamma\in\Gamma}\gamma Y,
\]
and that they are standard measure spaces endowed with the corresponding restrictions of $\sigma$, say $\mu=\sigma|_X$ and $\nu=\sigma|_Y$. 
The action of $\Gamma$ on $\Sigma/\Delta$ is isomorphic to the following action of $\Gamma$ on $X$: for each pair $(\gamma,x)\in\Gamma\times X$, there exists a unique element $\alpha(\gamma,x)\in\Delta$ such that $\gamma x\in X\alpha(\gamma,x)$. Hence the element $\gamma\cdot x:=\gamma x\alpha(\gamma,x)^{-1}$ belongs to $X$, and the mapping $(\gamma,x)\mapsto\gamma\cdot x$ defines an action of $\Gamma$ on $X$, and $\alpha$ is a Borel cocycle with values in $\Delta$:
\[
\alpha(\gamma_1\gamma_2,x)=\alpha(\gamma_1,\gamma_2\cdot x)\alpha(\gamma_2,x)\quad \forall x\in X,\ \forall\gamma_1,\gamma_2\in\Gamma.
\]

Similarly, the action of $\Delta$ on the orbit space $\Gamma\backslash\Sigma$
is isomorphic to the following action of $\Delta$ on $Y$: for each
pair $(y,\delta)\in Y\times\Delta$ there exists a unique
$\beta(y,\delta)\in\Gamma$ such that $y\delta\in
\beta(y,\delta)Y$. Then set
$y\cdot\delta=\beta(y,\delta)^{-1}y\delta\in Y$, so that this defines an
action of $\Delta$ on $Y$ on the right, and $\beta$ is a Borel
cocycle for this action, viz 
\[
\beta(y,\delta_{1})\beta(y\cdot\delta_{1},\delta_{2})
=\beta(y,\delta_{1}\delta_{2})\quad\forall y\in Y,\
\forall\delta_{1},\delta_{2}\in\Delta.
\]

Observe that $\mu$ (resp. $\nu$) is $\Gamma$-invariant (resp. $\Delta$-invariant), but that it need not be finite. However, there are invariant states on $L^{\infty}(X)$ and on $L^{\infty}(Y)$.

\par\vspace{3mm}
Theorem \ref{thmAME} will be a straightforward consequence of Theorem \ref{thmproper} for two reasons: the cocycle $\alpha$ is proper as the following lemma shows, and amenably measure equivalence is a symmetric property.

\begin{lemma}\label{3.1}
Retaining notation above, $\alpha$ is a proper cocycle from $\Gamma\times X$ to $\Delta$. More precisely:
\begin{enumerate}
\item [(a)] Let $1\in K\subset \Gamma$ be finite, let $A\subset X$ be a Borel set with finite measure and let $\varepsilon>0$. Then there exists a Borel set $A_{\varepsilon,K}\subset A$ and a finite set $F\subset \Delta$ such that $\sigma(A\smallsetminus A_{\varepsilon,K})<\varepsilon$ and
\[
\bigcup_{\gamma\in K}\gamma A_{\varepsilon,K}\subset \bigsqcup_{\delta\in F}X\delta.
\]
In particular $\alpha(K\times A_{\varepsilon,K})$ is a finite set and $\alpha$ satisfies condition (i) of Definition \ref{1.3}.
\item [(b)] Let $1\in F\subset\Delta$ be a finite set and let $\mathcal{A}_F$ be the set of Borel subsets $A$ of $X$ with finite measure for which there exists a finite set $K=K(A,F)\subset \Gamma$ such that 
\[
\bigsqcup_{\delta\in F}A\delta\subset \bigsqcup_{\gamma\in K}\gamma Y.
\]
Then the elements of $\mathcal{A}_F$ have the following properties:
\begin{enumerate}
\item [(b1)] For every Borel set $A\subset X$ with finite measure and for every $\varepsilon>0$, there exists $A_\varepsilon\in\mathcal{A}_F$ such that $A_\varepsilon\subset A$ and $\mu(A\smallsetminus A_\varepsilon)<\varepsilon$.
\item [(b2)] Let $A\in\mathcal{A}_F$ and let $K$ be a finite subset of $\Gamma$ such that 
\[
\bigsqcup_{\delta\in F}A\delta\subset \bigsqcup_{\gamma\in K}\gamma Y.
\]
For $\gamma\in \Gamma$, set, as in Definition \ref{1.3},
\[
X_{A}(\gamma)=\{x\in A\cap\gamma^{-1}\cdot A:\alpha(\gamma,x)\in F\}.
\]
Then $X_{A}(\gamma)=\emptyset$ for every $\gamma\notin KK^{-1}$.
\end{enumerate}
In particular, $\alpha$ satisfies condition (ii) of Definition \ref{1.3}.
\end{enumerate}
\end{lemma}
\textit{Proof.} (a) Using induction on $|K|$, it suffices to prove the claim for singleton sets. Thus, let $A\subset X$ be a Borel set with finite measure, let $\gamma\in \Gamma$ and let $\varepsilon>0$. There exists a finite set $F\subset\Delta$ such that 
\[
\sum_{\delta\notin F}\sigma(A\cap\gamma^{-1}X\delta)<\varepsilon.
\]
Then the Borel set $A_\varepsilon:=\bigsqcup_{\delta\in F}A\cap\gamma^{-1}X\delta$ is a subset of $A$ such that
$\sigma(A\smallsetminus A_\varepsilon)<\varepsilon$ and $\gamma A_\varepsilon\subset \bigsqcup_{\delta\in F}X\delta.$\\
(b1) Put $AF=\bigsqcup_{\delta\in F}A\delta$. Since $\sigma(AF)=|F|\sigma(A)=|F|\mu(A)<\infty$, there exists a finite set $K\subset\Gamma$ such that
\[
\sum_{\gamma\notin K}\sigma(AF\cap \gamma Y)<\varepsilon.
\]
Put $Z=AF\cap\left(\bigsqcup_{\gamma\in K}\gamma Y\right)$. For every $\delta\in F$, put $Z_\delta=(Z\delta^{-1})\cap X$, so that 
$
Z=\bigsqcup_{\delta\in F}Z_\delta\delta.
$
Finally, put 
\[
A_\varepsilon=\bigcap_{\delta\in F}Z_\delta.
\]
Then $A_\varepsilon\subset Z_\delta$ for every $\delta\in F$, hence $A_\varepsilon\delta\subset Z_\delta\delta\subset Z$ for every $\delta\in F$, so that
$\bigsqcup_{\delta\in F}A_\varepsilon\delta\subset Z$. One has:
\begin{eqnarray*}
\mu(A\smallsetminus A_\varepsilon) &=&
\sigma(A\smallsetminus A_\varepsilon)\\
&=&
\sigma\left(A\cap\left(\bigcup_{\delta\in F}Z_\delta^c\right)\right)\\
 &\leq &
 \sum_{\delta\in F}\sigma(A\cap Z_\delta^c)=\sum_{\delta\in F}\sigma(A\delta\smallsetminus Z_\delta\delta)\\
 &=&
 \sigma(AF\smallsetminus Z)<\varepsilon.
\end{eqnarray*}
This ends the first part of the proof of the lemma.\\
(b2) Let $A$ and $K$ be as stated, and let $\gamma\in\Gamma$. If $X_{A}(\gamma)$ contains some element $x$, then $\alpha(\gamma,x)\in F$ and $\gamma x\in AF$. It follows that 
\[
x\in\left(\bigsqcup_{\gamma'\in K}\gamma'Y\right)
\cap\left(\bigsqcup_{\gamma''\in K}\gamma^{-1}\gamma''Y\right).
\]
This implies that there are $\gamma',\gamma''\in K$ such that $\gamma=\gamma''\gamma'^{-1}\in KK^{-1}$.
\hfill $\square$

\bibliographystyle{amsplain}

\begin{thebibliography}{10}

\bibitem{AD}
{\sc C. Anantharaman-Delaroche.}
Amenable correspondences and approximation properties for von Neumann algebras.
{\em Pac. J. Math. 171\/} (1995), 309--341.

\bibitem{ccjjv}
{\sc P.-A. Cherix, M.~Cowling, P.~Jolissaint, P.~Julg, and A.~Valette.}
{\it Groups with the {H}aagerup property ({G}romov's
  a-{T}-menability).}
Birkh\"auser, {B}asel, (2001).

\bibitem{CH}
{\sc M. Cowling and U. Haagerup}.
Completely bounded multipliers of the Fourier algebra of a simple Lie group of real rank one.
{\em Invent. Math. 96\/} (1989), 507--549. 


\bibitem{Furman1}
{\sc A. Furman}
Gromov's measure equivalence and rigidity of higher rank lattices.
{\em Ann. Math. 150\/} (1999), 1059--1081.

\bibitem{Furman2}
{\sc A. Furman}
Orbit equivalence rigidity.
{\em Ann. Math. 150\/} (1999), 1083--1108.

\bibitem{Gre}
{\sc F.P. Greenleaf.}
Amenable actions of locally compact groups.
{\em Journal Funct. Anal. 4\/} (1969) 295--315.

\bibitem{Jolcoc}
{\sc P. Jolissaint.}
Borel cocycles, approximation properties and relative property T.
{\em Ergod. Th. \& Dynam. Syst. 20(2)\/} (2000) 483--499.

\bibitem{Jolinv}
{\sc P. Jolissaint.}
Invariant states and a conditional fixed point property for affine actions.
{\em Math. Ann. 304\/} (1996), 561--579.

\bibitem{Knudby}
{\sc S. Knudby.}
The weak Haagerup property.
{\em To appear in Trans. Amer. Math. Soc.\/} (2014), 41 pages.

\bibitem{Rudin}
{\sc W. Rudin.}
{\it Real and Complex Analysis.}
McGraw-Hill Inc. New-York, (1987).

\bibitem{Zimmer}
{\sc R.J. Zimmer.}
Amenable pairs of groups and ergodic actions and the associated von Neumann algebras.
{\em Trans. Amer. Math. Soc. 243\/} (1978), 271--286.

\end{thebibliography}

\end{document}